\documentclass[12pt]{amsart}

\usepackage[latin1]{inputenc}
\usepackage{amscd}
\usepackage{latexsym}
\usepackage{pstcol,pst-plot,pst-3d}
\usepackage{mathptmx}
\usepackage{multicol}
\psset{unit=0.7cm,linewidth=0.8pt,arrowsize=2.5pt 4}

\newpsstyle{fatline}{linewidth=1.5pt}
\newpsstyle{fyp}{fillstyle=solid,fillcolor=verylight}
\definecolor{verylight}{gray}{0.97}
\definecolor{light}{gray}{0.9}
\definecolor{medium}{gray}{0.85}
\def\frk{\frak}               

\def\Phi{{\frk n}}
\def\Phi{{\frk N}}
\def\xb{{\bold x}}
\def\yb{{\bold y}}
\def\zb{{\bold z}}
\def\opn#1#2{\def#1{\operatorname{#2}}} 
\opn\chara{char}
\opn\length{\ell}
\opn\pd{pd}
\opn\rk{rk}
\opn\projdim{proj\,dim}
\opn\injdim{inj\,dim}
\opn\rank{rank}
\opn\depth{depth}
\opn\grade{grade}
\opn\height{height}
\opn\embdim{emb\,dim}
\opn\codim{codim}

\opn\Tr{Tr}
\opn\bigrank{big\,rank}
\opn\superheight{superheight}\opn\lcm{lcm}
\opn\trdeg{tr\,deg}%
\opn\reg{reg}
\opn\lreg{lreg}
\opn\ini{in}
\opn\lpd{lpd}
\opn\div{div}
\opn\Div{Div}
\opn\cl{cl}
\opn\Cl{Cl}
\opn\Spec{Spec}
\opn\Supp{Supp}
\opn\supp{supp}
\opn\Sing{Sing}
\opn\Ass{Ass}
\opn\Ann{Ann}
\opn\Rad{Rad}
\opn\Soc{Soc}
\opn\Im{Im}
\opn\Ker{Ker}
\opn\Coker{Coker}
\opn\Am{Am}
\opn\Hom{Hom}
\opn\Tor{Tor}
\opn\Ext{Ext}
\opn\End{End}
\opn\Aut{Aut}
\opn\id{id}

\opn\nat{nat}
\opn\pff{pf}
\opn\Pf{Pf}
\opn\GL{GL}
\opn\SL{SL}
\opn\mod{mod}
\opn\ord{ord}
\opn\Gin{Gin}
\opn\aff{aff}
\opn\con{conv}
\opn\relint{relint}
\opn\st{st}
\opn\lk{lk}
\opn\cn{cn}
\opn\core{core}
\opn\vol{vol}
\opn\link{link}
\opn\star{star}
\opn\com{com}
\opn\gr{gr}

\def\pot#1#2{#1[\kern-0.28ex[#2]\kern-0.28ex]}

\opn\dirlim{\underrightarrow{\lim}}
\opn\inivlim{\underleftarrow{\lim}}
\let\union=\cup
\let\sect=\cap
\let\dirsum=\oplus

\let\iso=\cong

\let\Sect=\bigcap
\let\Dirsum=\bigoplus

\let\to=\rightarrow
\let\To=\longrightarrow
\def\Implies{\ifmmode\Longrightarrow \else
        \unskip${}\Longrightarrow{}$\ignorespaces\fi}
\def\implies{\ifmmode\Rightarrow \else
        \unskip${}\Rightarrow{}$\ignorespaces\fi}
\def\iff{\ifmmode\Longleftrightarrow \else
        \unskip${}\Longleftrightarrow{}$\ignorespaces\fi}

\let\:=\colon
\newtheorem{Theorem}{Theorem}[section]
\newtheorem{Lemma}[Theorem]{Lemma}
\newtheorem{Corollary}[Theorem]{Corollary}
\newtheorem{Proposition}[Theorem]{Proposition}

\let\epsilon\varepsilon
\let\phi=\varphi
\let\kappa=\varkappa
\def\qed{\ifhmode\textqed\fi
      \ifmmode\ifinner\quad\qedsymbol\else\dispqed\fi\fi}
\def\textqed{\unskip\nobreak\penalty50
       \hskip2em\hbox{}\nobreak\hfil\qedsymbol
       \parfillskip=0pt \finalhyphendemerits=0}
\def\dispqed{\rlap{\qquad\qedsymbol}}

\opn\dis{dis}
\def\pnt{{\raise0.5mm\hbox{\large\bf.}}}

\opn\Lex{Lex}

\begin{document}

\title{Distributive Lattices, Bipartite Graphs and Alexander Duality}
\author{J\"urgen Herzog and  Takayuki Hibi}
\address{J\"urgen Herzog, Fachbereich Mathematik und
Informatik, Universit\"at Duisburg-Essen, Campus Essen,
45117 Essen, Germany}
\email{juergen.herzog@uni-essen.de}
\address{Takayuki Hibi, Department of Pure and Applied Mathematics,
Graduate School of Information Science and Technology,
Osaka University, Toyonaka, Osaka 560-0043, Japan}
\email{hibi@math.sci.osaka-u.ac.jp}
\maketitle
\begin{abstract}
A certain squarefree monomial ideal $H_P$ arising from
a finite partially ordered set $P$ will be studied
from viewpoints of both commutative algebra and combinatorics.
First, it is proved that the defining ideal
of the Rees algebra of $H_P$ possesses
a quadratic Gr\"obner basis.  Thus in particular all powers
of $H_P$ have linear resolutions.  Second, the minimal free graded
resolution of $H_P$ will be constructed explicitly and
a combinatorial formula to compute the Betti numbers of $H_P$
will be presented.  Third, by using the fact that
the Alexander dual of the simplicial complex $\Delta$
whose Stanley--Reisner ideal coincides with $H_P$
is Cohen--Macaulay, all the Cohen--Macaulay bipartite graphs
will be classified.
\end{abstract}

\section*{Introduction}
Let $P$ be a finite partially ordered set ({\em poset} for short)
and write ${\mathcal J}(P)$ for the finite poset which consists of all
poset ideals of $P$, ordered by inclusion.
Here a {\em poset ideal} of $P$ is
a subset $I$ of $P$ such that if $x \in I$, $y \in P$
and $y \leq x$, then $y \in I$.
In particular
the empty set as well as $P$ itself is a poset ideal of $P$.
It follows easily that ${\mathcal J}(P)$ is a finite distributive
lattice
\cite[p.\ 106]{StanleyEnumerative}.
%
%
Conversely, Birkhoff's fundamental structure theorem
\cite[Theorem 3.4.1]{StanleyEnumerative}
guarantees that, for any finite distributive lattice ${\mathcal L}$,
there exists a unique poset $P$ such that ${\mathcal L} = {\mathcal
J}(P)$.

Let $P$ be a finite poset with $|P| = n$,
where $|P|$ is the cardinality of $P$,
and let $S = K[ \{ x_p, y_p \}_{p \in P}]$ denote
the polynomial ring in $2n$ variables over a field $K$
with each $\deg x_p = \deg y_p = 1$.

We associate each poset ideal $I$ of $P$
with the squarefree monomial
\[
u_I
= (\prod_{p \in I} x_p) (\prod_{p \in P \setminus I} y_p)
\]
of $S$ of degree $n$.
In particular $u_P = \prod_{p \in P} x_p$ and
$u_\emptyset = \prod_{p \in P} y_p$.
 
The normal affine semigroup ring
$K[\{ u_I \}_{I \in {\mathcal J}(P)}]$
is studied in \cite{HibiDistributiveLattice} from viewpoints of
%
%
both commutative algebra and combinatorics.

In the present paper, however, we are interested in
the squarefree monomial ideal
\[
H_P = ( \{ u_I \}_{I \in {\mathcal J}(P)})
\]
of $S$
generated by all $u_I$ with $I \in {\mathcal J}(P)$.

The outline of the present paper is as follows.
First, in Section $1$ we study the Rees algebra
${\mathcal R}(H_P)$ of $H_P$ and establish our fundamental
Theorem \ref{quadraticbinomials} which says that
the defining ideal of ${\mathcal R}(H_P)$ possesses
a reduced Gr\"obner basis consisting
of quadratic binomials whose initial monomials are squarefree.
Thus ${\mathcal R}(H_P)$ turns out to be
normal and Koszul (Corollary \ref{normalandKoszul}),
and
all powers of $H_P$ have linear resolutions
(Corollary \ref{allpowershavelinearresolutions}).

Second, in Section $2$
the minimal graded free $S$-resolution of $H_P$
is constructed explicitly.
See Theorem \ref{resolution}.
The resolution tells us how to compute the Betti numbers
$\beta_i(H_P)$ of $H_P$
in terms of the combinatorics of the distributive lattice
${\mathcal L} = {\mathcal J}(P)$.
In fact, if $b_i({\mathcal L})$ is the number of intervals
$[I, J]$ of
${\mathcal L} = {\mathcal J}(P)$ which are Boolean lattices
of rank $i$, then the $i$th Betti number
$\beta_i(H_P)$ of $H_P$ coincides with $b_i({\mathcal L})$.
See Corollary \ref{interpretation}.
(A Boolean lattice of rank $i$ is the
distributive lattice $B_i$ which consists of all subsets
of $\{1, \ldots, i \}$, ordered by inclusion.)
Thus in particular for a finite distributive lattice
${\mathcal L} = {\mathcal J}(P)$, one has
$\sum_{i\geq 0}(-1)^ib_i({\mathcal L})=1$.
See Corollary \ref{alternating}.
In addition, it is shown that the ideal $H_P$ is of height 2 and
a formula to compute the multiplicity of $S/H_P$ will be given.
See Proposition \ref{multiplicity}
(and Corollary \ref{formula}).

Let $\Delta_P$ denote the simplicial complex
on the vertex set
$\{ x_p, y_p \}_{p \in P}$
such that the squarefree monomial ideal $H_P$
coincides with the Stanley--Reisner ideal $I_{\Delta_P}$.
In Section $3$ the Alexander dual $\Delta_P^\vee$ of $\Delta_P$
will be studied.
Since the Stanley--Reisner ideal $H_P = I_{\Delta_P}$
has a linear resolution, it follows from
\cite[Theorem 3]{EagonReiner}
%
%
that $\Delta_P^\vee$ is Cohen--Macaulay.
It will turn out that the Stanley--Reisner ideal
$I_{\Delta_P^\vee}$ of $\Delta_P^\vee$ is an edge ideal
of a finite bipartite graph.
Somewhat surprisingly,
this simple observation enables us to classify all Cohen--Macaulay
bipartite graphs.  In fact, Theorem \ref{classification}
says that a finite bipartite graph $G$ is Cohen--Macaulay
if and only if $G$ comes from the comparability graph
of a finite poset.

\section{Monomial ideals arising from distributive lattices}
Work with the same notation as in Introduction.
Let $P$ be a finite poset with $|P| = n$
and $S = K[\{ x_p, y_p \}_{p \in P}]$
the polynomial ring in $2n$ variables over a field $K$
with each $\deg x_p = \deg y_p = 1$.
Recall that we associate each poset ideal $I$ of $P$
with the squarefree monomial $u_I
= (\prod_{p \in I} x_p) (\prod_{p \in P \setminus I} y_p)$
of $S$ of degree $n$,
and introduce the ideal
$H_P = ( \{ u_I \}_{I \in {\mathcal J}(P)})$
of $S$.

Let ${\mathcal R}(H_P)$ denote the Rees algebra of $H_P$
and ${\mathcal W}_P$ the defining ideal of ${\mathcal R}(H_P)$.
In other words, ${\mathcal R}(H_P)$ is the affine semigroup ring
\[
{\mathcal R}(H_P) = K[\{ x_p, y_p \}_{p \in P},
\{ u_I t\}_{I \in {\mathcal J}(P)}]
\, \, \, \, \, ( \subset
K[\{ x_p, y_p \}_{p \in P}, t] )
\]
and ${\mathcal W}_P$ is the kernel of the surjective ring homomorphism
$\varphi : K[\xb, \yb, \zb] \to {\mathcal R}(H_P)$, where
\[
K[\xb, \yb, \zb] = K[\{ x_p, y_p \}_{p \in P},
\{ z_I \}_{I \in {\mathcal J}(P)}]
\]
is the polynomial ring over $K$ and where $\varphi$
is defined by setting
$\varphi(x_p) = x_p$,
$\varphi(y_p) = y_p$ and
$\varphi(z_I) = u_I t$.

For the convenience of our discussion, in the remainder of
the present section, we will use the notation
$P = \{ p_1, \ldots, p_n \}$ and write $x_i$, $y_i$
instead of $x_{p_i}$, $y_{p_i}$.
Let $<_{lex}$ denote the lexicographic order on $S$ induced by
the ordering $x_1 > \cdots > x_n > y_1 > \cdots >y_n$ and
$<^{\sharp}$ the reverse lexicographic order on
$K[\{ z_I \}_{I \in {\mathcal J}(P)}]$ induced by an ordering of
the variables $z_I$'s such that
$z_I > z_J$ if $J < I$ in ${\mathcal J}(P)$.
We then introduce the new monomial order $<_{lex}^{\sharp}$ on $T$
by setting
\[
(\prod_{i=1}^{n} x_i^{a_i} y_i^{b_i})
(z_{I_1} \cdots z_{I_q})
<_{lex}^{\sharp}
(\prod_{i=1}^{n} x_i^{a'_i} y_i^{b'_i})
(z_{I'_1} \cdots z_{I'_{q'}})
\]
if either

(i) \ \
$\prod_{i=1}^{n} x_i^{a_i} y_i^{b_i}
<_{lex}
\prod_{i=1}^{n} x_i^{a'_i} y_i^{b'_i}$

\noindent
or

(ii) \ \
$\prod_{i=1}^{n} x_i^{a_i} y_i^{b_i}
=
\prod_{i=1}^{n} x_i^{a'_i} y_i^{b'_i}$
and
$z_{I_1} \cdots z_{I_q}
<^{\sharp}
z_{I'_1} \cdots z_{I'_{q'}}$.

\begin{Theorem}
\label{quadraticbinomials}
The reduced Gr\"obner basis ${\mathcal G}_{<_{lex}^{\sharp}}({\mathcal
W}_P)$
of the defining ideal
${\mathcal W}_P \subset K[\xb, \yb, \zb]$
with respect to the monomial order $<_{lex}^{\sharp}$ consists
of quadratic binomials whose initial monomials are squarefree.
\end{Theorem}

\begin{proof}
The reduced Gr\"obner basis of
${\mathcal W}_P \cap K[\{ z_I \}_{I \in {\mathcal J}(P)}]$
with respect to the reverse lexicographic order $<^{\sharp}$
coincides with
${\mathcal G}_{<_{lex}^{\sharp}}({\mathcal W}_P)
\cap K[\{ z_I \}_{I \in {\mathcal J}(P)}]$.
It follows from \cite{HibiDistributiveLattice}
that ${\mathcal G}_{<_{lex}^{\sharp}}({\mathcal W}_P)
\cap K[\{ z_I \}_{I \in {\mathcal J}(P)}]$
consists of those binomials
\[
z_I z_J - z_{I \wedge J} z_{I \vee J}
\]
such that $I$ and $J$ are incomparable in the distributive
lattice ${\mathcal J}(P)$.

It is known \cite[Corollary 4.4]{Sturmfels} that
%
%
the reduced Gr\"obner basis of ${\mathcal W}_P$ consists of irreducible
binomials of $K[\xb, \yb, \zb]$.
Let
\[
f = (\prod_{i=1}^{n} x_i^{a_i} y_i^{b_i})
(z_{I_1} \cdots z_{I_q}) -
(\prod_{i=1}^{n} x_i^{a'_i} y_i^{b'_i})
(z_{I'_1} \cdots z_{I'_{q}})
\]
be an irreducible binomial of
$K[\xb, \yb, \zb]$
belonging to
${\mathcal G}_{<_{lex}^{\sharp}}({\mathcal W}_P)$
with
\[
(\prod_{i=1}^{n} x_i^{a_i} y_i^{b_i})
(z_{I_1} \cdots z_{I_q})
\]
its initial monomial, where
$I_1 \leq \cdots \leq I_q$ and
$I'_1 \leq \cdots \leq I'_{q}$.

Let $f \not\in K[\{ z_I \}_{I \in {\mathcal J}(P)}]$.
Let $j$ denote an integer for which
$I'_j \not\subset I_j$.
Such an integer exists.
In fact, if $I'_j \subset I_j$ for all $j$,
then each $a_i = 0$ and each $b'_i = 0$.
This is impossible since
$(\prod_{i=1}^{n} x_i^{a_i} y_i^{b_i})
(z_{I_1} \cdots z_{I_q})$ is the initial monomial
of $f$.

Let $p_i \in I'_{j} \setminus I_j$.
Then $p_i$ belongs to each of
$I'_j, I'_{j+1}, \ldots, I'_q$, and
does not belong to
each of $I_1, I_2, \ldots, I_j$.
Hence $a_i > 0$.

Let $p_{i_0} \in P$ with
$p_{i_0} \in I'_{j} \setminus I_j$
such that
$I_j \cup \{ p_{i_0} \} \in {\mathcal J}(P)$.
Thus $a_{i_0} > 0$.
Let $J = I_j \cup \{ p_{i_0} \}$.
Then the binomial $g = x_{i_0} z_{I_j} - y_{i_0} z_J$
belongs to ${\mathcal W}_P$ with $x_{i_0} z_{I_j}$ its initial
monomial.
Since $x_{i_0} z_{I_j}$ divides the initial monomial of $f$,
it follows that the initial monomial of $f$ must coincides with
$x_{i_0} z_I$, as desired.
\end{proof}

It is well known that a homogeneous affine semigroup ring
whose defining ideal has an initial ideal which is
generated by squarefree
(resp.\ quadratic) monomials is normal (resp.\ Koszul).
See, e.g., \cite[Proposition 13.15]{Sturmfels} and
\cite{KoszulAlgebra}.
%
%

\begin{Corollary}
\label{normalandKoszul}
Let $P$ be an arbitrary finite poset.
Then the Rees algebra ${\mathcal R}(H_P)$ is normal and Koszul.
\end{Corollary}

On the other hand, Stefan Blum \cite{JuergenStudent}
%
%
proved that if the Rees algebra of an ideal is Koszul, then
all powers of the ideal have linear resolutions.

\begin{Corollary}
\label{allpowershavelinearresolutions}
Let $P$ be an arbitrary finite poset.
Then all powers of $H_P$ have linear resolutions.
\end{Corollary}

\section{The free resolution and Betti numbers of $H_P$}

Corollary \ref{allpowershavelinearresolutions} says that
the monomial ideal $H_P$ arising from
a finite poset $P$ has a linear resolution.
The main purpose of the present section is to
construct a minimal graded free $S$-resolution ${\Bbb F}={\Bbb F}_P$
of $H_P$ explicitly.

Let $P$ be a finite poset with $|P| = n$
and $S = K[\{ x_p, y_p \}_{p \in P}]$
the polynomial ring in $2n$ variables over a field $K$
with each $\deg x_p = \deg y_p = 1$.
Recall that, for each poset ideal $I$ of $P$,
we associate the squarefree monomial
$u_I = (\prod_{p \in I} x_p) (\prod_{p \in P \setminus I} y_p)$
of $S$ of degree $n$.
Let $H_P$ denote the ideal of $S$
generated by all $u_I$
with $I \in {\mathcal J}(P)$.

The maximal elements of a poset ideal $I$ of $P$ are called
the {\em generators} of $I$.
Let $M(I)$ denote the set of generators of $I$.

The construction of
a minimal graded free $S$-resolution ${\Bbb F}={\Bbb F}_P$ of $H_P$
is achieved as follows:
For all $i\geq 0$ let ${\Bbb F}_i$ denote the free $S$-module with
basis
\[
e(I,T),
\]
where
\[
I\in {\mathcal J}(P),\ T\subset P,\
I\sect T\subset M(I),\ |I\sect T|=i
\, \, \,
\text{and}
\, \, \,
|I\union T|=n+i.
\]
Extending the partial order on $P$ to a total order,
we define for $i>0$  the differential
\[
\partial:\ {\Bbb F}_{i}\to {\Bbb F}_{i-1}
\]
  by
\[
\partial(e(I,T))=\sum_{p\in I\sect T}
(-1)^{\sigma(I\sect T,p)}(x_{p}e(I\setminus\{p\},T)-
y_{p}e(I,T\setminus\{p\})),
\]
where for a subset $Q\subset P$ and $p\in Q$
we set $\sigma(Q,p)=|\{q\in Q\: q<p\}|$.

With the notation introduced we have

\begin{Theorem}
\label{resolution}
The complex ${\Bbb F}$ is a graded minimal free $S$-resolution of
$H_P$.
\end{Theorem}

\begin{proof}
We define an augmentation $\epsilon\: {\Bbb F}_0\to H_P$ by setting
\[
\epsilon(e(I,T))=u_I
\]
for all $e(I,T)\in {\Bbb F}_0$. Note that if $e(I,T)$ is
a basis element of ${\Bbb F}_0$, then
$T=[n]\setminus I$, so that $\epsilon$  is well defined.

We first show that
\[
\begin{CD}
\cdots  @>\partial >>
{\Bbb F}_1 @>\partial >> {\Bbb F}_0 @>\epsilon >> H_P @>>> 0
\end{CD}
\]
is a complex.

Let $e(I,T)\in {\Bbb F}_1$ with $I\sect T=\{p\}$. Then
\begin{eqnarray*}
(\epsilon\circ \partial)(e(I,T))
& = & x_p\epsilon(e(I\setminus\{p\}, T))
-y_p\epsilon(e(I,T\setminus \{p\})) \\
& = & x_pu_{I\setminus\{p\}}-y_pu_{I}=0.
\end{eqnarray*}
Thus $\partial\circ\epsilon=0$, as desired.

Next we show that $\partial\circ\partial=0$.
Let $e(I,T)\in {\Bbb F}_{i+1}$
and set $L=I\sect T$.
Then
\begin{eqnarray*}
\partial\circ\partial(e(I,T))&=&
\sum_{p\in
L}(-1)^{\sigma(L,p)}\bigl(x_p\partial(e(I\setminus\{p\},T))\\
& & \, \, \, \, \,
- y_p\partial(e(I,T\setminus\{p\})\bigr)\\
&=&\sum_{p\in L}(-1)^{\sigma(L,p)}\bigl[x_p\bigl(\sum_{q\in L,q\neq p}
(-1)^{\sigma(L\setminus\{p\},q)}\bigl(x_qe(I\setminus\{p,q\},T)\\
& & \, \, \, \, \, -y_qe(I\setminus\{p\},T\setminus\{q\}))\bigr)\\
& & \, \, \, \, \,
-y_p\bigl(\sum_{q\in L, q\neq p}(-1)^{\sigma(L\setminus
\{p\},q)}(x_qe(I\setminus\{q\},T\setminus\{p\})\\
& & \, \, \, \, \,
-y_qe(I,T\setminus\{p,q\}))\bigr)\bigr]\\
&=&\sum_{p,q\in L, p\neq q}(-1)^{\sigma(L,p)+\sigma(L\setminus\{p\},q)}
x_px_qe(I\setminus\{p,q\},T)\\
& & \, \, \, \, \,
-\sum_{p,q\in L, p\neq  q}(-1)^{\sigma(L,p)+\sigma(L\setminus\{p\},q)}
x_py_qe(I\setminus\{p\},T\setminus\{q\})\\
& & \, \, \, \, \,
-\sum_{p,q\in L, p\neq  q}
(-1)^{\sigma(L,p)\sigma(L\setminus\{p\},q)}x_qy_pe(I\setminus\{q\},
T\setminus\{p\})\\
& & \, \, \, \, \,
+\sum_{p\in L,p\neq q}(-1)^{\sigma(L,p)+\sigma(L\setminus\{p\},q)}
y_py_qe(I,T\setminus\{p,q\})\\
&=&0.
\end{eqnarray*}
The last equality holds since
$(-1)^{\sigma(L,p)+\sigma(L\setminus\{p\},q)}=-(-1)^{\sigma(L,q)+
\sigma(L\setminus\{q\},p)}$.

In order to prove that the augmented complex
\[
\begin{CD}
  \cdots @>>> {\Bbb F}_1@>\partial >> {\Bbb F}_0@> \epsilon >> H_P@>>>
0
\end{CD}
\]
is exact we show:
\begin{enumerate}
\item[(1)] $H_0({\Bbb F})=H_P$,
\item[(2)] ${\Bbb F}$ is acyclic.
\end{enumerate}

For the proof of (1) we note that the Taylor relations
\[
r_{I,J}=x_{J\setminus I}y_{I\setminus J}e(I)-x_{I\setminus J}
y_{J\setminus I}e(J),\qquad I,J\in {\mathcal J}(P)
\]
generate the first syzygy module of $H_P$.
Here we set for simplicity $e(I)$ for the basis element
$e(I,P\setminus I)$ in ${\Bbb F}_0$, and denote by
$x_Ay_B$ the monomial $\prod_{p\in A}x_p\prod_{q\in B}y_q$.

Observe that
\[
r_{I,J}=x_{J\setminus I}r_{I,I\sect J}-x_{I\setminus J}r_{J,I\sect J}.
\]
Hence it suffices to show that $r_{I,J}\in\partial({\Bbb F}_1)$
for all $I,J\in L$  with $J\subset
I$. To this end we choose a sequence
$J=I_0\subset I_1\subset\ldots I_{m-1}\subset I_m=I$ of poset
ideals such that $I_j=I_{j-1}\union \{p_j\}$ for $j=1,\dots,m$.
Then
\[
r_{I,J}=\sum_{j=1}^m(\prod_{k=j+1}^mx_{p_k}\prod_{k=1}^{j-1}y_{p_k})
r_{I_j,I_{j-1}}.
\]
The assertion follows since $r_{I_j,I_{j-1}}=-\partial(e(I_j,P\setminus
I_{j-1}))$ for all $j$.

\medskip

We prove (2), that is, the acyclicity of $\Bbb F$ by induction on
$|P|$.
If $P=\{p\}$, then
$H_P=(x_p,y_p)$, and ${\Bbb F}$ can be identified with the Koszul
complex
associated with
$\{x_p,y_p\}$, and hence is acyclic.

Suppose now that $|P|>1$. Let $q\in P$ be a maximal element and
let $Q$ be the subposet
$P\setminus\{q\}$.
We define a map
\[
\phi\: {\Bbb F}_Q\to {\Bbb F}_P,\qquad e_i(I,T)\mapsto
e_i(I,T\union\{q\})
\]
It is clear that $\phi$ is an injective map of complexes whose induced
map
$H_Q=H_0({\Bbb F}_Q)\to
H_0({\Bbb F}_P)=H_P$ is multiplication by $y_q$.
Let ${\Bbb G}$ be the quotient complex ${\Bbb
F}_P/{\Bbb F}_Q$. Since the multiplication map is injective,
the short exact sequence of complexes
\[
\begin{CD}
0@>>> {\Bbb F}_Q@>>> {\Bbb F}_P@>>> {\Bbb G}@>>> 0
\end{CD}
\]
induces the long exact homology sequence
\[
\begin{CD}
\cdots @>>> H_2({\Bbb G})@>>>
H_1({\Bbb F}_Q)@>>> H_1({\Bbb F}_P)@>>> H_1({\Bbb G})@>>> 0
\end{CD}
\]
By induction hypothesis, $H_i({\Bbb F}_Q)=0$ for $i>0$.
Hence it suffices
to show that $H_i({\Bbb G})=0$ for $i>0$.

The principal order ideal $(q)$ consists of all $p\in P$
with $p\leq p$.
Let $R$ be the subposet
$P\setminus (q)$, and let ${\Bbb C}$ be the mapping cone of the complex
homomorphism
\[
\begin{CD}
{\Bbb F}_R@> -y_q >> {\Bbb F}_R.
\end{CD}
\]
Then we get an exact sequence
\[
\begin{CD}
0@>>> {\Bbb F}_R@>>> {\Bbb C}@>>> {\Bbb F}_R[-1]@>>> 0
\end{CD}
\]
Here ${\Bbb F}_R[-1]$ is the complex ${\Bbb F}_R$ shifted to
the `left', that is,
$({ \Bbb
F}_R[-1])_i= ({\Bbb F}_R)_{i-1}$ for all $i$.
 
By our induction hypothesis ${\Bbb F}_R$ is  acyclic.
Thus from the long exact sequence
\[
\begin{CD}
H_1({\Bbb C})@>>>H_0({\Bbb F}_R) @>-y_q>>
H_0({\Bbb F}_R)@>>> H_0({\Bbb C})@>>> 0\\
\cdots @>>> H_2({\Bbb C})@>>>H_1({\Bbb F}_R) @> -y_q>>
H_1({\Bbb F}_R)@>>>
\end{CD}
\]
we deduce that $H_i({\Bbb C})=0$ for $i>1$.
We also get $H_1({\Bbb C})=0$,
since $H_0({\Bbb F}_R)=H_R$, and since multiplication by $y_q$ is
injective on $H_R$.
Thus we see that ${\Bbb C}$  is acyclic.

We now claim that ${\Bbb C}\iso {\Bbb G}$, thereby proving that
${\Bbb G}$ is acyclic,
as desired.

In order to prove this claim we first notice that
${\Bbb C}_i=({\Bbb F}_R)_{i-1}\dirsum ({\Bbb
F}_R)_i$  for $i\geq 0$ (where $({\Bbb F}_R)_{-1}=0$).
Thus if $r=|R|$, then
${\Bbb C}_i$ has the
basis ${\mathcal C}_i = {\mathcal B}_{i-1}\union {\mathcal B}_i$, where
\begin{eqnarray*}
{\mathcal B}_i&=&\{e(I,T)\: I\in L(R),\ T\subset R,\
I\sect T\subset M(I),
|I\sect T|=i,\\ &&\ |I\union T|=r+i\}.
\end{eqnarray*}
On the other hand ${\Bbb G}_i$ has the basis
\begin{eqnarray*}
{\mathcal G}_i&=&\{e(I,T)\: I\in L(P),\ (q)\subset I,\
T\subset P,\ I\sect T\subset M(I),
|I\sect T|=i,\\ &&\ |I\union T|=n+i\}.
\end{eqnarray*}

Let $\psi_i: {\Bbb C}_i \to {\Bbb G}_i$ be the $S$-linear homomorphism
with
\[
\psi_i(e(I,T))=\left\{
\begin{array}{ll}
e(I\union (q),T\union\{q\}) &\text{if $e(I,T)\in {\mathcal
B}_{i-1}$};\\
e(I\union (q),T)& \text{if $e(I,T)\in {\mathcal B}_{i}$}.\\
\end{array}
\right.
\]
It is easy to see that all  $\psi_i$ are bijections and induce
an isomorphism of complexes.
\end{proof}

Suppose $P$ is of cardinality $n$ and $P$ is an antichain,
i.e., any two elements of
$P$ are incomparable.  Then $B_n={\mathcal J}(P)$ is
called the {\em Boolean lattice of rank $n$}.

Let now ${\mathcal L}$ be an arbitrary finite distributive lattice,
and let $I,J\in {\mathcal L}$ with $I\leq J$. Then the
set
\[
[I,J]=\{M\in {\mathcal L}\: I\leq M\leq J\}
\]
is called an {\em interval} in ${\mathcal L}$. The interval $[I,J]$
with the induced partial order is again a
distributive lattice.
Let $b_i({\mathcal L})$ denote the number of intervals of
${\mathcal L}$ which are isomorphic to
Boolean lattices of rank $i$.
In particular, $b_0({\mathcal L}) = |{\mathcal L}|$.
These numbers have an algebraic interpretation.

Recall that for a
graded $S$-module $M$,
\[
\beta_{i}(M)=\dim_K  \Tor_i^S(M,K)
\]
is called the {\em $i$th Betti-number of $M$}.
If ${\Bbb F}$ is a graded minimal free resolution of $M$,
then $\beta_i(M)$ is nothing but the rank of ${\Bbb F}_i$.

\begin{Corollary}
\label{interpretation}
Let $P$ be a finite poset, ${\mathcal L} = {\mathcal J}(P)$ the
distributive lattice
and $H_P$ the squarefree monomial ideal arising from $P$.
Then
\begin{enumerate}
\item[(a)] $b_i({\mathcal L})=\beta_i(H_P)$ for all $i$;
\item[(b)] the following three numbers are  equal:
\begin{enumerate}
\item[(i)] the projective dimension of $H_P$;
\item[(ii)] the maximum of the ranks of Boolean lattices which are
isomorphic
to an interval of ${\mathcal L}$;
\item[(iii)] the Sperner number of $P$, i.e., the maximum of the
cardinalities of
antichains of $P$.
\end{enumerate}
\end{enumerate}
\end{Corollary}

\begin{proof}
(a) For each $i\geq 0$, let ${\mathcal S}_i$ be the set of pairs
$(I,S)$, where $I\in{\mathcal L}$,  $S\subset M(I)$ and $|S|=i$, and let ${\mathcal B}_i$ be the set of
basis elements  $e(I,T)$ of ${\Bbb F}_i$. Then
\[
{\mathcal B}_i\To {\mathcal S}_i,\quad e(I,T)\mapsto (I,I\sect T)
\]
establishes a bijection between these two sets.

Since for each $(I,S)\in {\mathcal S}_i$, the elements in $S$ are
pairwise incomparable it is clear
that $[I\setminus S,I]$ is isomorphic to a Boolean lattice
of rank $i$.

Conversely,  suppose $[J,I]$
is isomorphic to a Boolean lattice of rank $i$. Then $S=I\setminus J$
is of a set of cardinality $i$, and $J\union T\in{\mathcal L}$ for all subsets $T\subset S$.

Suppose that $S\not\subset M(I)$. Then there exists, $q\in S$ and $p\in
I$ such that $p>q$. If $p\in J$, then $q\in J$, a contradiction. Thus $p\in S$, and hence $(J,p)\in
{\mathcal L}$. This is again a contradiction, because it would imply that $q\in (J,p)$. Hence we have
shown that $(I,S)\in {\mathcal S}_i$.

It follows that the assignment $e(I,T)\mapsto [I\setminus(I\sect T),I]$
establishes a bijection between the basis of ${\Bbb F}_i$ and the intervals of $[J,I]$ in
${\mathcal L}$ which are isomorphic to Boolean lattices.

(b) is an immediate consequence of (a) and its proof.
\end{proof}

\begin{Corollary}
\label{alternating}
Let ${\mathcal L}$ be a finite distributive lattice.
Then
\[
\sum_{i\geq 0}(-1)^{i+1}b_i({\mathcal L})=1.
\]
\end{Corollary}

Let $\Delta_P$ be the simplicial complex attached to the squarefree
monomial
ideal $H_P$. In the next section we will see
(Lemma \ref{flag})
that the Stanley--Reisner ideal
attached to the Alexander dual $\Delta_P^{\vee}$ is
generated by the monomials $x_py_q$ such that $p\leq q$.
Hence for the Stanley--Reisner ideal
of $\Delta_P$ we have
\[
I_{\Delta_P}=\Sect_{p,q\in P,\ p\leq q}(x_p,x_q).
\]
In particular we get
\begin{Proposition}
\label{multiplicity}
Let $P$ be a finite poset.
Then the squarefree monomial ideal $H_P$ is of height 2, and the
multiplicity of
$S/H_P$ is given by
\[
e(S/H_P)=|\{(p,q)\: p,q\in P,\ p\leq q\}|.
\]
\end{Proposition}

Let $I\subset S$ be an arbitrary graded ideal with graded minimal free
resolution
\[
\begin{CD}
0@>>> \Dirsum_{j=1}^{\beta_s}S(-a_{sj})@>>> \cdots @>>>
\Dirsum_{j=1}^{\beta_1}S(-a_{1j})@>>> S@>>>
S/I@>>> 0.
\end{CD}
\]
Suppose the height of $I$  equals $h$.
Then by a formula of Peskine and Szpiro \cite{PS} one has
\[
e(S/I)=\frac{(-1)^h}{h!}\sum_{i=1}^s(-1)^i\sum_{j=1}^{\beta_i}a_{ij}^h.
\]
Applying this formula in our situation and using Corollary
\ref{interpretation}
and Proposition \ref{multiplicity} we get

\begin{Corollary}
\label{formula}
Let $P$ be a finite poset with $|P|=n$, and let ${\mathcal L} =
{\mathcal J}(P)$ be
the distributive lattice.  Then
\[
|\{(p,q)\: p,q\in P,\ p\leq q\}|
=\frac{1}{2}\sum_{i\geq 0}(-1)^{i+1}b_{i}({\mathcal L})(n+i)^2.
\]
\end{Corollary}
 
We close this section with an example.
Let $P$ be the poset with Hasse diagram

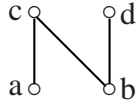
\begin{center}
\psset{unit=0.5cm}
\begin{pspicture}(0,0)(3,3)
  \psline(0.5,0.6)(0.5,2.35)
  \psline(0.55,2.4)(2.4,0.55)
  \psline(2.5,0.6)(2.5,2.35)
  \rput(0.5,2.5){$\circ$}
  \rput(0.5,0.45){$\circ$}
  \rput(2.5,0.45){$\circ$}
  \rput(2.5,2.5){$\circ$}
  \rput(0,0.5){a}
  \rput(3,0.5){b}
  \rput(0,2.5){c}
  \rput(3,2.5){d}
\end{pspicture}
\end{center}

\noindent
The distributive lattice ${\mathcal L} = {\mathcal J}(P)$ has the Hasse
diagram

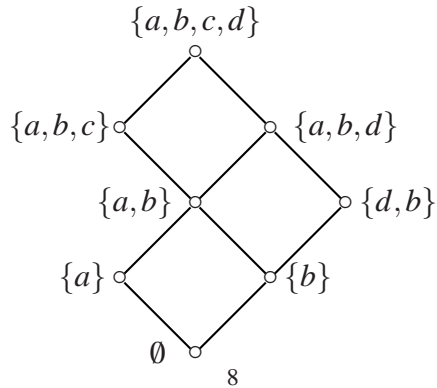
\begin{center}
\psset{unit=1cm}
\begin{pspicture}(-1,0)(4,5)
  \psline(1.06,0.06)(1.94,0.94)
  \psline(2.06,1.06)(2.94,1.94)
  \psline(0.94,0.06)(0.06,0.94)
  \psline(0.06,1.06)(0.94,1.94)
  \psline(1.94,1.06)(1.06,1.94)
  \psline(1.06,2.06)(1.94,2.94)
  \psline(2.94,2.06)(2.06,2.94)
  \psline(0.94,2.06)(0.06,2.94)
  \psline(1.94,3.06)(1.06,3.94)
  \psline(0.06,3.06)(0.94,3.94)
  \rput(1,0){$\circ$}
  \rput(2,1){$\circ$}
  \rput(3,2){$\circ$}
\rput(0,1){$\circ$}
\rput(1,2){$\circ$}
\rput(2,3){$\circ$}
\rput(0,3){$\circ$}
\rput(1,4){$\circ$}
\rput(0.5,0){$\emptyset$}
\rput(-0.5,1){$\{a\}$}
\rput(0.2,2){$\{a,b\}$}
\rput(-0.8,3){$\{a,b,c\}$}
\rput(1,4.4){$\{a,b,c,d\}$}
\rput(3,3){$\{a,b,d\}$}
\rput(3.7,2){$\{d,b\}$}
\rput(2.5,1){$\{b\}$}
\end{pspicture}
\end{center}

\noindent
Thus $H_P=(uvwx,avwx,buwx,abwx,bduw,abcx,abdw,abcd)$.
Here we use for
convenience the indeterminates $a,b,c,d,u,v,w,x$ instead of $x_p$ and
$y_p$.
The free resolution of $H_P$ is given by
\[
\begin{CD}
0@>>> S^3(-6)@>>> S^{10}(-5)@>>> S^8(-4)@>>> H_P@>>> 0.
\end{CD}
\]
We see from the Hasse diagram that the $i$th Betti number of $H_P$
coincides with number of intervals
of ${\mathcal L}$ which are isomorphic to Boolean lattices of rank $i$.
The number of pairs $(p,q)$ in the poset $P$ with $p\leq q$
is equal to 7, and this is also the
number we get from Corollary \ref{formula},
namely $(1/2)(-8\cdot 16+10\cdot 25-3\cdot 36)=7$.

\section{Alexander duality and Cohen--Macaulay bipartite graphs}
We refer the reader to, e.g., \cite{StanleyGreenBook},
\cite{BrunsHerzog}
and \cite{HibiAlgebraicCombinatorics}
%
%
for fundamental information about Stanley--Reisner rings.

Let $P = \{ p_1, \ldots, p_n \}$ be a finite poset
and
$S = K[x_1, \ldots, x_n, y_1, \ldots, y_n]$
the
polynomial ring in $2n$ variables over a field
$K$ with each $\deg x_i = \deg y_i = 1$.
We will use the notation $x_i$, $y_i$ instead of
$x_{p_i}$, $y_{p_i}$, and set $V_n = \{ x_1, \ldots, x_n, y_1, \ldots, y_n \}$.

Recall that $H_P$ is the ideal of $S$ which is
generated by those squarefree monomials
$u_I = (\prod_{p_i \in I} x_i)(\prod_{p_i \in P \setminus I} y_i)$
with $I \in {\mathcal J}(P)$.  It then follows that
there is a unique simplicial complex $\Delta_P$
on $V_n$ such that
the {\em Stanley--Reisner ideal} $I_{\Delta_P}$
coincides with $H_P$.
We study the {\em Alexander dual} $\Delta_P^\vee$
of
$\Delta_P$, which is the simplicial complex
\[
\Delta_P^\vee = \{ V_n \setminus F \, : \, F \not\in \Delta_P \}
\]
on $V_n$.

\begin{Lemma}
\label{flag}
The Stanley--Reisner ideal of $\Delta_P^\vee$ is generated by
those squarefree quadratic monomials $x_i y_j$ such that
$p_i \leq p_j$ in $P$.
\end{Lemma}

\begin{proof}
Let $w = x_1 \cdots x_n y_1 \cdots y_n$. If $u$ is a squarefree
monomial of $S$, then we write $\supp(u)$ for the support of $u$,
i.e., $\supp(u)
= \{ x_i \, : \, x_i \, \, \mbox{divides} \, \, u \} \cup
\{ y_j \, : \, y_j \, \, \mbox{divides} \, \, u \}$. Now
since $\{ \supp(u_I) \, : \, I \in {\mathcal J}(P) \}$
is the set of minimal nonfaces of $\Delta_P$, it follows that
$\{ \supp(w / u_I) \, : \, I \in {\mathcal J}(P) \}$
is the set of facets (maximal faces) of $\Delta_P^\vee$.
Our work is to find the minimal nonfaces of $\Delta_P^\vee$.
Since $\supp(w / u_\emptyset) = x_1 \cdots x_n$ and
$\supp(w / u_P) = y_1 \cdots y_n$,
both $\{ x_1, \ldots, x_n \}$ and $\{ y_1, \ldots, y_n \}$
are faces of $\Delta_P^\vee$.
Let $F \subset V_n$ be a nonfaces of $\Delta_P^\vee$.
Let $F_x = F \cap \{ x_1, \ldots, x_n \}$ and
$F_y = \{ x_j \, : \, y_j \in F \}$.
Then $F_x \neq \emptyset$ and $F_y \neq \emptyset$.
Since $\{ x_i, y_i \}$ is a minimal nonface of $\Delta_P^\vee$,
we will assume that $F_x \cap F_y = \emptyset$.
Since $F$ is a nonface, there exists {\em no} poset ideal
$I$ of $P$ with $F_x \cap \{ x_i \, : \, p_i \in I \}
= \emptyset$
and $F_y \subset \{ x_i \, : \, p_i \in I \}$.
Hence there are $x_i \in F_x$
and $x_j \in F_y$ such that $p_i < p_j$.
Thus $\{ x_i, y_j \}$ is a nonface of $\Delta_P^\vee$.
Hence the set of minimal nonfaces of $\Delta_P^\vee$
consists of those $2$-element subsets $\{ x_i, y_j \}$ of $V_n$
such that $p_i \leq p_j$ in $P$, as required.
\end{proof}

Let $G$ be a finite graph on the vertex set
$[N] = \{ 1, \ldots, N \}$
with no loops and no multiple edges.
We will assume that $G$ possesses no isolated vertex, i.e.,
for each vertex $i$ there is an edge $e$ of $G$ with $i \in e$.
A {\em vertex cover} of $G$ is a subset $C \subset [N]$
such that, for each edge $\{ i, j \}$ of $G$,
one has either $i \in C$ or $j \in C$.
Such a vertex cover $C$ is called {\em minimal}
if no subset $C' \subset C$ with $C' \neq C$ is a vertex cover
of $G$.
We say that a finite graph $G$ is {\em unmixed} if all minimal
vertex covers of $G$ have the same cardinality.

Let $K[\zb] = K[z_1, \ldots, z_N]$ denote the polynomial ring
in $N$ variables over a field $K$.
The {\em edge ideal} of $G$ is the ideal $I(G)$ of $K[\zb]$
generated by those squarefree quadratic monomials $z_i z_j$
such that $\{ i, j \}$ is an edge of $G$.
A finite graph $G$ on $[N]$ is called {\em Cohen--Macaulay} over $K$
if the quotient ring $K[\zb]/I(G)$ is Cohen--Macaulay.
Every Cohen--Macaulay graph is unmixed
(\cite[Proposition 6.1.21]{MonomialAlgebras}).
%
%

A finite graph $G$ on $[N]$ is {\em bipartite} if
there is a partition $[N] = W \cup W'$ such that each edge
of $G$ is of the form $\{ i , j \}$ with $i \in W$ and $j \in W'$.
A basic fact on the graph theory says that a finite graph
$G$ is bipartite if and only if $G$ possesses no cycle of odd length.
A {\em tree} is a connected graph with no cycle.
A tree is Cohen--Macaulay if and only if it is unmixed
(\cite[Corollary 6.3.5]{MonomialAlgebras}).

Given a finite poset $P = \{ p_1, \ldots, p_n \}$,
we write $G(P)$ for the bipartite graph
on the vertex set
$\{ x_1, \ldots, x_n \} \cup \{ y_1, \ldots, y_n \}$
whose edges are those $\{ x_i, y_j \}$ such that
$p_i \leq p_j$ in $P$.
Lemma \ref{flag} says that
the Stanley--Reisner ideal of
$\Delta_P^\vee$ is equal to the edge ideal of $G(P)$.
Since the Stanley--Reisner ideal $H_P = I_{\Delta_P}$
has a linear resolution, it follows from
\cite[Theorem 3]{EagonReiner}
%
%
that $\Delta_P^\vee$ is Cohen--Macaulay.
Then
\cite[Theorem 6.4.7]{MonomialAlgebras}
says that $\Delta_P^\vee$ is shellable.
Hence $I_{\Delta_P}$
has linear quotients (e.g., \cite{HHZAlexanderDuality}).
%
%

\begin{Corollary}
\label{linearquotients}
The ideal $H_P$ has linear quotients.
\end{Corollary}

We now turn to the problem of classifying the
Cohen--Macaulay bipartite graphs
by using the Alexander dual
$\Delta_P^\vee$.

Let $G$ be a finite bipartite graph on the vertex set
$W \cup W'$
with $W = \{ i_1, \ldots, i_s \}$
and $W' = \{ j_1, \ldots, j_t \}$,
where $s \leq t$.
For each subset $U$ of $W$, we write $N(U)$ for the set of
those vertices $j \in W'$
for which there is a vertex $i \in U$
such that $\{ i, j \}$ is an edge of $G$.
The well-known ``marriage theorem'' in graph theory
says that if $|U| \leq |N(U)|$
for all subsets $U$ of $W$, then
there is a subset
$W'' = \{ j_{\ell_1}, \ldots, j_{\ell_s} \} \subset
W'$ with $|W''| = s$ such that
$\{ i_k, j_{\ell_k} \}$ is an edge of $G$ for $k = 1, 2, \ldots, s$.

Let $G$ be a finite bipartite graph on the vertex set
$W \cup W'$ and suppose that $G$ is unmixed.
Since each of $W$ and $W'$
is a minimal vertex cover, one has $|W| = |W'|$.
Let $W = \{ x_1, \ldots, x_n \}$ and
$W' = \{ y_1, \ldots, y_n \}$.
Since
$(W \setminus U) \cup N(U)$
is a vertex cover of $G$
for all subsets $U$ of $W$
and
since $G$ is unmixed, it follows that
$|U| \leq |N(U)|$
for all subsets $U$ of $W$.
Thus the marriage theorem enables us to
assume that $G$ satisfies the condition
as follows:
$(\sharp)$
$\{ x_i, y_i \}$
is an edge of $G$ for all $1 \leq i \leq n$.

\begin{Lemma}
\label{wonderful}
Work with the same notation as above and, furthermore,
suppose that $G$ is a Cohen--Macaulay graph. Then,
after a suitable change of the labeling of variables
$y_1, \ldots, y_n$, the edge set of $G$
satisfies the condition
$(\sharp)$ together with the condition
as follows:
$(\sharp \sharp)$ if $\{ x_i, y_j \}$
is an edge of $G$, then $i \leq j$.
\end{Lemma}

\begin{proof}
Let $\Delta$ be the Cohen--Macaulay complex on the vertex set
$W \cup W'$ whose Stanley--Reisner ideal $I_{\Delta}$
coincides with $I(G)$.
Recall that every Cohen--Macaulay complex is strongly connected
and that all links of a Cohen--Macaulay complex are
again Cohen--Macaulay.
Since both $W$ and $W'$ are facets of $\Delta$,
it follows (say, by induction on $n$) that,
after a suitable change of the labeling of variables
$x_1, \ldots, x_n$ and
$y_1, \ldots, y_n$, the subset
$F_i = \{ y_1, \ldots, y_i, x_{i+1}, \ldots, x_n \}$
is a facet of $\Delta$ for each $0 \leq i \leq n$,
where $F_0 = W$ and $F_n = W'$.
In particular $\{ x_i, y_j \}$
cannot be an edge of $G$ if $j < i$.
In other words, the edge set of $G$ satisfies the conditions
$(\sharp)$ and $(\sharp \sharp)$, as required.
\end{proof}

\begin{Theorem}
\label{classification}
Let $G$ be a finite bipartite graph on the vertex set $W \cup W'$,
where $W = \{ x_1, \ldots, x_n \}$ and
$W' = \{ y_1, \ldots, y_n \}$, and suppose that the edge set of $G$
satisfies the conditions $(\sharp)$ and $(\sharp \sharp)$.
Then $G$ is a Cohen--Macaulay graph if and only if the following
condition $(\sharp \sharp \sharp)$ is satisfied:

\medskip

$(\sharp \sharp \sharp)$ If $\{ x_i, y_j \}$ and $\{ x_j, y_k \}$
are edges of $G$ with $i < j < k$, then
$\{ x_i, y_k \}$ is an edge of $G$.
\end{Theorem}

\begin{proof}
{\bf (``Only if'')}
Let $G$ be a Cohen--Macaulay graph
satisfying $(\sharp)$ and $(\sharp \sharp)$ and
$\Delta$ the Cohen--Macaulay complex on the vertex set
$W \cup W'$ whose Stanley--Reisner ideal coincides with $I(G)$.
Let $\{ x_i, y_j \}$ and $\{ x_j, y_k \}$
be edges of $G$ with $i < j < k$ and suppose that
$\{ x_i, y_k \}$ is not an edge of $G$.
Since every Cohen--Macaulay complex is pure and since
$\{ x_i, y_k \}$ is a face of $\Delta$, it follows that
there is an $n$-element subset $F \subset W \cup W'$
of $G$ with
$\{ x_i, y_k \} \subset F$ such that $F$ is independent in $G$,
i.e., no $2$-element subset of $F$ is an edge of $G$.
One has $y_j \not\in F$ and $x_j \not\in F$
since $\{ x_i, y_j \}$ and $\{ x_j, y_k \}$ are edges of $G$.
Since $\{ x_\ell, y_\ell \}$ is an edge of $G$ for each
$1 \leq \ell \leq n$, the independent subset $F$ can contain
both $x_i$ and $y_i$ for no $1 \leq i \leq n$.
Thus to find such an $n$-element independent set $F$ is
impossible.

\medskip

\noindent
{\bf (``If'')} Now, suppose that a finite bipartite graph
$G$ on the vertex set $W \cup W'$
satisfies the conditions $(\sharp)$, $(\sharp \sharp)$
together with $(\sharp \sharp \sharp)$.
Let $\leq$ denote the binary relation on
$P = \{ p_1, \ldots, p_n \}$ defined by setting
$p_i \leq p_j$ if $\{ x_i, y_j \}$ is an edge of $G$.
By $(\sharp)$ one has
$p_i \leq p_i$ for each $1 \leq i \leq n$.
By $(\sharp \sharp)$ if $p_i \leq p_j$ and $p_j \leq p_i$,
then $p_i = p_j$.
By $(\sharp \sharp \sharp)$
if $p_i \leq p_j$ and $p_j \leq p_k$, then
$p_i \leq p_k$.
Thus $\leq$ is a partial order on $P$.
Lemma \ref{flag} then guarantees that $G = G(P)$.
Hence $G$ is Cohen--Macaulay, as desired.
\end{proof}
 
\begin{Corollary}
\label{pureandstronglyconnected}
Let $G$ be a finite bipartite graph and $\Delta$ the simplicial
complex whose Stanley--Reisner ring coincides with $I(G)$.
Then $G$ is Cohen--Macaulay
if and only if $\Delta$ is pure and strongly connected.
\end{Corollary}

Work with the same situation as in the ``if'' part of
the proof of Theorem \ref{classification}.
Let $\com(P)$ denote the {\em comparability graph} of $P$,
i.e., $\com(P)$ is the finite graph on $\{ p_1, \ldots, p_n \}$
whose edges are those $\{ p_i, p_j \}$ with $i \neq j$
such that $p_i$ and $p_j$ are comparable in $P$.
It then follows from
\cite[pp. 184 -- 185]{MonomialAlgebras}
that the Cohen--Macaulay type of the Cohen--Macaulay ring
$S / I(G)$, where $S = K[x_1, \ldots, x_n, y_1, \ldots, y_n]$,
is the number of maximal independent subsets of $\com(P)$,
i.e., the number of maximal antichains of $P$.
Hence $G$ is Gorenstein, i.e.,
$S / I(G)$ is a Gorenstein ring, if and only if $P$ is
an antichain.

\begin{Corollary}
\label{Gorensteinbipartitegraph}
A Cohen-Macaulay bipartite graph $G$ is Gorenstein
if and only if $G$ is the disjoint union of edges.
\end{Corollary}


\begin{thebibliography}{99}

\bibitem{JuergenStudent}
S.\ Blum, Subalgebras of bigraded Koszul algebras,
{\em J. Algebra} {\bf 242} (2001), 795 -- 809.

\bibitem{BrunsHerzog}
W.\ Bruns and J.\ Herzog, ``Cohen--Macaulay rings,''
Revised Edition, Cambridge University Press, 1996.

\bibitem{EagonReiner}
J.\ Eagon and V.\ Reiner, Resolutions of
Stanley--Reisner rings and Alexander duality,
{\em J. Pure Appl. Algebra} {\bf 130} (1998),
265--275.

\bibitem{KoszulAlgebra}
R.\ Fr\"oberg, Koszul algebras,
{\em in} ``Advances in Commutative Ring Theory''
(D.\ E.\ Dobbs, M.\ Fontana and S.-E.\ Kabbaj, Eds.),
Lecture Notes in Pure and Appl. Math.,
Volume 205, Dekker, New York, NY, 1999,
pp. 337 -- 350.

\bibitem{HHZAlexanderDuality}
J.\ Herzog, T.\ Hibi and X.\ Zheng,
Dirac's theorem on chordal graphs and Alexander duality,
preprint, 2003.

\bibitem{HibiDistributiveLattice}
T.\ Hibi, Distributive lattices, affine semigroup
rings and algebras with straightening laws,
{\em in} ``Commutative Algebra and Combinatorics''
(M.\ Nagata and H.\ Matsumura, Eds.),
Advanced Studies in Pure Math.,
Volume 11, North--Holland, Amsterdam, 1987,
pp. 93 -- 109.

\bibitem{HibiAlgebraicCombinatorics}
T.\ Hibi, ``Algebraic Combinatorics on Convex Polytopes,''
Carslaw, Glebe, N.S.W., Australia, 1992.
 
\bibitem{PS}
C.\ Peskine and L.\ Szpiro, Syzygies and multiplicities,
{\em C.\ R.\ Acad.\ Sci.\ Paris.\ S\'er.\ A} {\bf 278} (1974),
1421 -- 1424.

\bibitem{StanleyEnumerative}
R.\ P.\ Stanley, ``Enumerative Combinatorics, Volume I,''
Wadsworth \& Brooks/Cole, Monterey, CA, 1986.

\bibitem{StanleyGreenBook}
R.\ P.\ Stanley, ``Combinatorics and Commutative Algebra,''
Second Edition, Birkh\"auser, Boston, MA, 1996.

\bibitem{Sturmfels}
B.\ Sturmfels, ``Gr\"obner Bases and Convex Polytopes,''
Amer. Math. Soc., Providence, RI, 1995.

\bibitem{MonomialAlgebras}
R.\ H.\  Villareal, ``Monomial Algebras,''
Dekker, New York, NY, 2001.

\end{thebibliography}
\end{document}